\documentclass[11pt]{amsart}
\usepackage{amscd,amssymb}
\usepackage{amsthm,amsmath,amssymb}
\usepackage[matrix,arrow]{xy}

\newtheorem{theorem}[equation]{Theorem}

\newtheorem{proposition}[equation]{Proposition}
\newtheorem{lemma}[equation]{Lemma}
\newtheorem{corollary}[equation]{Corollary}
\newtheorem{conjecture}[equation]{Conjecture}

\theoremstyle{definition}
\newtheorem{example}[equation]{Example}

\theoremstyle{remark}
\newtheorem{remark}[equation]{Remark}

 \makeatletter\@addtoreset{equation}{section}
 \makeatother

\newcommand{\Supp}{\operatorname{Supp}}

\newcommand{\Sing}{\operatorname{Sing}}

\newcommand{\mult}{\operatorname{mult}}

\newcommand{\rank}{\operatorname{rank}}

\newcommand{\ZZ}{{\mathbb Z}}

\newcommand{\QQ}{{\mathbb Q}}
\newcommand{\PP}{{\mathbb P}}

\newcommand{\MMM}{\mathcal M}
\newcommand{\OOO}{\mathcal O}
\newcommand{\DDD}{\mathcal D}

\newcommand{\III}{\mathcal I}
\newcommand{\LLL}{\mathcal L}

\begin{document}

\Large
\begin{center}
\textbf{Factorial hypersurfaces in $\PP^4$ with nodes}
\end{center}
\vspace{5mm}

\normalsize
\begin{center}
\textbf{Ivan Cheltsov$^*$ and Jihun Park$^{**}$
}\\
\vspace{5mm} \small{\itshape $^*$School of Mathematics, University
of
Edinburgh,  Mayfield Road, Edinburgh EH9 3JZ, UK.}\\
\vspace{5mm} \small{\itshape $^{**}$Department of Mathematics,
POSTECH, Pohang, Kyungbuk, $790$-$784$, Republic of Korea.}
\end{center}
\vspace{5mm}

\hrulefill \vspace{2mm}

 \small {\textbf{Abstract.} We prove that for $n=5, 6, 7$ a nodal hypersurface of degree $n$
in $\PP^4$ is factorial if it has at most $(n-1)^2-1$ nodes.
 \vspace{2mm}

  }

\hrulefill \normalsize
 \footnotetext[1]{Corresponding author : J. Park, Tel.
 +82-54-279-2059 Fax. +82-54-279-2799.}
 \footnotetext[2]{E-mail addresses :
\texttt{cheltsov@yahoo.com} (I. Cheltsov);
\texttt{wlog@postech.ac.kr} (J. Park)}

\section{Introduction}
\label{section:introduction} Unless otherwise mentioned, every
variety is always assumed to be projective, normal, and defined
over $\mathbb{C}$. We also consider every divisor in the linear
system $|\mathcal{O}_{\mathbb{P}^{n}}(k)|$ as a hypersurface in
$\mathbb{P}^{n}$ for simplicity.

A variety $X$ is called factorial if each Weil divisor of $X$ is
Cartier. The factoriality is a very subtle property. It depends on
both the local types of singularities and their global position.
Also, it depends on the field of definition of the variety. In the
present paper, we study the factoriality of a hypersurface in
$\PP^4$. However, we confine our consideration to the case when
they have only simple double points, \emph{i.e.}, nodes.

Let $V_n$ be  a nodal hypersurface of degree $n$ in $\PP^4$. Then
the Picard group is isomorphic to the 2nd integral cohomology
because $H^1(W, \OOO_W)=H^2(W, \OOO_W)=0$ on a resolution $W$ of
$V_n$.  In this case,  the variety $V_n$ is factorial if and only
if the global topological property
\[\rank (H^2(V_n, \ZZ))=\rank (H_4(V_n, \ZZ)),\]
holds. Note that the duality mentioned above fails on singular
varieties in general. The nodes on $V_n$ may have an effect on the
integral (co)homology groups of $V_n$ (See \cite{Cy01}). However,
the rank of the 2nd integral cohomology group of $V_n$ is 1  by
Lefschetz theorem. Therefore, to determine whether the 3-fold
hypersurface $V_n$ is factorial or not, we have to see whether the
rank of the 4th integral homology group of $V_n$ is $1$ or not.
But, it is not simple to compute the rank of the 4th integral
homology group of $V_n$. Fortunately, the paper \cite{Cy01} gives
us a great method to compute the rank of the 4th integral homology
group of $V_n$, which reduce the topological problem to a rather
simple combinatorial problem. To be precise, the rank of the 4th
integral homology group  of $V_n$ can be obtained by the following
way:

\begin{theorem}\label{theorem-of-Cynk} Let $V_n$ be a nodal hypersurface of degree $n$ in $\PP^4$.
The rank of the $4$th integral homology group $H_4(V_n, \ZZ)$ is
equal to
\[\#|\Sing (V_n)| - I+1,\]
where $I$ is the number of independent conditions which vanishing
on $\Sing (V_n)$ imposes on homogeneous forms of degree $2n-5$ on
$\PP^4$.
\end{theorem}
\begin{proof} See \cite{Cy01}. \end{proof}

Therefore, the hypersurface $V_n$ is factorial if and only if the
set of nodes of the hypersurface $V_n$  is
$(2n-5)$-normal\footnote{In general, a subscheme $X$ of $\PP^N$ is
called $d$-normal in $\PP^N$ if the first cohomology of the sheaf
of ideal of $X$ twisted by $\OOO_{\PP^N}(d)$ is zero. Throughout
this paper, we consider a finite set of points in $\PP^n$ as a
zero-dimensional  reduced subscheme of $\PP^n$.} in $\PP^4$ (see
\cite{Cy01}), in other words, the singular points of the
hypersurface $V_n$ impose linearly in\-de\-pen\-dent conditions on
hypersurfaces  of degree $2n-5$ in $\mathbb{P}^{4}$.

The geometry of the hypersurface $V_n$ crucially depends on its
factoriality. For example, in the case $n=4$ the hypersurface
$V_n$ is  non-rational whenever it is factorial (see \cite{Me03}),
which is not true without the factoriality condition.

Let us show an easy way to get a non-factorial hypersurface.

\begin{example}
\label{example:non-factorial-quintic} Suppose that $V_n$ is given
by the equation
$$
x_0g(x_0,x_1,x_2,x_3,x_4)+x_1f(x_0,x_1,x_2,x_3,x_4)=0\subset\mathbb{P}^{4}\cong\mathrm{Proj}(\mathbb{C}[x_0,x_1,x_2,x_3,x_4]),%
$$
where $g$ and $f$ are general homogeneous polynomials of degree
$n-1$. Then the hypersurface $V_n$ has exactly $(n-1)^2$ nodes.
They are located on the $2$-plane defined by $x_0=x_1=0$. The
hypersurface $V_n$ is not factorial because the hyperplane section
$x_0=0$ splits into two non-Cartier divisors while the Picard
group is generated by a hyperplane section (see \cite{AnFra59}).
\end{example}

On the other hand, in the case when $\mathrm{Sing}(V_n)<(n-1)^2$,
every smooth surface in $V_n$ is cut by a hypersurface in
$\mathbb{P}^{4}$ due to \cite{CiGe03}. Therefore, it is natural
that we should expect the following to be true.

\begin{conjecture}
\label{conjecture:factoriality} Every nodal hypersurface of degree
$n$ in $\PP^4$ with at most $(n-1)^2-1$ nodes is factorial.
\end{conjecture}

 Con\-jec\-ture~\ref{conjecture:factoriality} for $n=2$ and $3$ is trivial. For the
case $n=4$ Conjecture~\ref{conjecture:factoriality} is proved in
\cite{Ch04e}.  In this paper we prove the following result.

\begin{theorem}
\label{theorem:main}%
Conjecture~\ref{conjecture:factoriality} holds for $n=5, 6$, and
$7$.
\end{theorem}

Note that the following result is proved in \cite{Ch04t}.

\begin{theorem}
\label{theorem:Cheltsov}%
A nodal hypersurface of degree $n$ in $\PP^4$ with at most
$\frac{(n-1)^2}{4}$ nodes is factorial.
\end{theorem}

Therefore, at least asymptotically
Conjecture~\ref{conjecture:factoriality} is not far from being
true. To our surprise, the conjecture below implies
Conjecture~\ref{conjecture:factoriality}.

\begin{conjecture}
\label{conjecture:projections} Let $\Sigma$ be a subset in
$\mathbb{P}^{4}$ such that at most $k(n-1)$ points in $\Sigma$ can
be contained in a curve of degree $k$, where $|\Sigma|<(n-1)^{2}$.
Then at most $k(n-1)$ points in $\phi_4(\Sigma)$ can be contained
in a curve of degree $k$ in $\mathbb{P}^{2}$, where
$\phi_4:\mathbb{P}^{4}\dasharrow\mathbb{P}^{2}$ is a general
projection.
\end{conjecture}

Unfortunately, we are unable to prove
Conjecture~\ref{conjecture:factoriality}
 now, but we believe
that the proof of Theorem~\ref{theorem:main}
 can help us to find new approaches to a
proof of Conjecture~\ref{conjecture:factoriality}.

\bigskip

\emph{Acknowledgements.} The authors would like to thank
I.\,Aliev, A.\,Corti, M.\,Gri\-nen\-ko, V.\,Iskov\-s\-kikh, K. Oh,
Yu.\,Pro\-kho\-rov, and V.\,Sho\-ku\-rov  for invaluable
conversations. The second author has been supported by KRF Grant
2005-070-C00005 in Republic of Korea.

\section{Preliminaries}
\label{section:preliminaries}

\subsection{Projections and linear systems with zero-dimensional base loci}

Let $X$ be a smooth variety with a $\QQ$-divisor
$B_{X}=\sum_{i=1}^{k}a_{i}B_{i}$, where $B_{i}$ is a prime divisor
on $X$ and $a_{i}$ is a positive rational number. Let $\pi:Y\to X$
be a birational morphism of a smooth variety $Y$ such that the
union of all the proper transforms of the divisors $B_{i}$ and all
the $\pi$-exceptional divisors forms a divisor with simple normal
crossing on $Y$.

Let $B_{Y}$ be the proper transform of $B_{X}$ on $Y$. Put
$B^{Y}=B_{Y}-\sum_{i=1}^{n}d_{i}E_{i}$, where each $E_{i}$ is an
exceptional divisor of the morphism $\pi$ and $d_{i}$ is the
rational number such that the equivalence
$$
K_{Y}+B^{Y}\sim_{\mathbb{Q}}\pi^{*}\big(K_{X}+B_{X}\big)
$$
holds. Then the log pair $(Y, B^{Y})$ is called the log pull back
of the log pair $(X, B_{X})$ with respect to $\pi$, while the
number $d_{i}$ is called the discrepancy of the log pair $(X,
B_{X})$ with respect to the $\pi$-exceptional divisor $E_i$.

By $\mathcal{L}(X, B_{X})$ we denote the subscheme of the variety
$X$ associated to the ideal sheaf
$$
\mathcal{I}\big(X, B_{X}\big)=\pi_{*}\Big(\mathcal{O}_{Y}\big(\lceil -B^{Y}\rceil\big)\Big),%
$$
which  is called the log canonical singularity subscheme of the
log pair $(X, B_{X})$.

We then obtain the following result due to \cite{Sh92}.

\begin{theorem}
\label{theorem:Shokurov} Suppose that the divisor $K_{X}+B_{X}+H$
is numerically equivalent to a Cartier divisor for some  nef and
big $\mathbb{Q}$-divisor $H$ on the variety $X$. Then for every
$i>0$ we have
$$
H^{i}\Big(X, \mathcal{I}(X, B_{X})\otimes
\mathcal{O}_{X}(K_{X}+B_{X}+H)\Big)=0.
$$
\end{theorem}

Theorem~\ref{theorem:Shokurov} gives us  a useful tool to study
the normality of a finite set in $\PP^N$.

\begin{lemma}\label{lemma:non-vanishing}
Let $\MMM$ be a linear system consisting of hypersurfaces of
degree $k$ on $\PP^N$. If the base locus $\Lambda$ of the linear
system $\MMM$ is zero-dimensional, then the finite set $\Lambda$
is $N(k-1)$-normal in $\PP^N$.
\end{lemma}

\begin{proof}  Let $H_1, \cdots, H_r$ be general divisors in the linear
system $\MMM$, where $r$ is sufficiently big. We put
$$B=\frac{N}{r}\sum_{i=1}^r H_i.$$

Then the log pair $(\PP^N, B)$ is klt in the outside of the base
locus $\Lambda$. For each point $p\in \Lambda$, we have
$\mult_pB\geq N$. Therefore, $\Supp(\LLL(\PP^N, B))=\Lambda$.

Since the divisor $K_{\PP^N}+B+H$ is $\QQ$-linearly equivalent to
$N(k-1)H$, where $H$ is a hyperplane,
we obtain $H^1(\PP^N, \III(\PP^N,
B)\otimes\OOO_{\PP^N}(N(k-1)))=0$ from
Theorem~\ref{theorem:Shokurov}. Because $\Supp(\LLL(\PP^N,
B))=\Lambda$ and the scheme $\LLL(\PP^N, B)$ is zero-dimensional,
the set $\Lambda$ that is the reduced scheme of $\LLL(\PP^N, B)$
must be $N(k-1)$-normal in $\PP^N$.
\end{proof}

Let $\Sigma$ be a finite set of points in $\PP^N$, $N\geq 3$, such
that no $k(d-1)+1$ points of $\Sigma$ lie on a curve of degree $k$
in $\PP^N$ for each $k\geq 1$, where $d\geq 3$ is a fixed integer.
Fix a $2$-plane $\Pi$ in $\PP^N$. We consider the projection
\[\phi_N:\PP^N\dasharrow\Pi\cong\PP^2\]
from a general $(N-3)$-dimensional linear space $L$ onto the
$2$-plane $\Pi$.
\begin{lemma}\label{lemma:zero-dimensional}
Let $\Lambda$ be a subset of $\Sigma$ and let $\MMM$ be the linear
system of hypersurfaces in $\PP^N$ of degree $k$ that contains
$\Lambda$. If $|\Lambda|>k(d-1)$ but the set $\phi_N(\Lambda)$ is
contained in an irreducible curve on $\Pi$ of degree $k$, then the
base locus of the linear system $\MMM$ is zero-dimensional.
\end{lemma}

\begin{proof} Suppose that the base locus of $\MMM$ contains an
irreducible curve $Z$. Let $C$ be an irreducible curve on $\Pi$ of
degree $k$ that contains $\phi_N(\Lambda)$ and let $W$ be the cone
in $\PP^N$ over the curve $C$ with vertex $L$. Since $W$ is a
hypersurface of degree $k$ in $\PP^N$ containing the set
$\Lambda$, it belongs to the linear system $\MMM$. In particular,
the curve $Z$ is contained in the hypersurface $W$. Therefore, the
curve $Z$ is mapped onto the curve $C$ because the linear space
$L$ is general and the curve $C$ is irreducible. The curve $Z$ has
degree $k$ because the restriction $\phi_N|_Z$ is a birational
morphism to $C$.

If there is a point $p$ in $\Lambda\setminus Z$, then the
projection $\phi_N$ maps the point $p$  to the outside of $C$
because of the generality of the projection $\phi_N$. Therefore,
the set $\Lambda$ must be contained in $Z$ because
$\phi_N(\Lambda)$ is contained in $C$. However, the curve $Z$
cannot contain more than $k(d-1)$ points of $\Sigma$.
\end{proof}

\begin{corollary}\label{corollary:line} A line on $\Pi$ contains at most $d-1$ points of
$\phi_N(\Sigma)$.
\end{corollary}
\begin{proof} It immediately follows from
Lemma~\ref{lemma:zero-dimensional}. \end{proof}

\begin{corollary}\label{corollary:conic} For $N=3$, a curve of degree $k$ on $\Pi$ contains at most $k(d-1)$ points of
$\phi_3(\Sigma)$ if $d\geq k^2+1$.
\end{corollary}

\begin{proof} For $k=1$, it is true because of
Corollary~\ref{corollary:line}. Assume that the claim is true for
$k<\ell$. We then suppose that there is a subset $\Lambda$ of
$\Sigma$ such that $|\Lambda|>\ell (d-1)$ and the image
$\phi_3(\Lambda)$ lie on a  curve $C$ of degree $\ell$ on $\Pi$.
The curve $C$ must be irreducible because of our assumption.
Therefore, it follows from Lemma~\ref{lemma:zero-dimensional} that
the base locus of the linear system $\MMM$ of  hypersurfaces of
degree $\ell$ in $\PP^3$ containing the set $\Lambda$ is
zero-dimensional. Let $Q_1$, $Q_2$, and $Q_3$ be general members
in $\MMM$. Then we obtain a contradictory inequality
\[\ell^3=Q_1\cdot Q_2\cdot Q_3> \ell (d-1)\geq \ell^3.\]
Therefore, for $d\geq k^2+1$, a curve of degree $k$ on $\Pi$
contains at most $k(d-1)$ points of $\phi_3(\Sigma)$.
\end{proof}

\begin{corollary}\label{corollary:conic2}
For $N=4$, a curve of degree $k$ on $\Pi$ contains at most
$k(d-1)$ points of $\phi_4(\Sigma)$ if $d\geq k^2+1$.
\end{corollary}
\begin{proof}
Let $\alpha :\PP^4\dasharrow\PP^3$ be the projection from a
generic point $o_1\in\PP^4$. We first claim that for $d\geq k^2+1$
no $k(d-1)+1$ points of $\alpha(\Sigma)$ lie on a  curve of degree
$k$ in $\PP^3$. It is obviously true for $k=1$. Assume that the
claim is true for $k<\ell$.  And then we suppose that there is a
subset $\Lambda$ of $\ell(d-1)+1$ points in $\Sigma$ such that
$\alpha(\Lambda)$ lie on a  curve $C$ of degree $\ell$ in $\PP^3$.
The curve $C$ must be irreducible. Let $\MMM$ be the linear system
of  hypersurfaces of degree $\ell$ in $\PP^4$ passing through all
the points of $\Lambda$. Then the proof of
Lemma~\ref{lemma:zero-dimensional} shows the base locus is
zero-dimensional. Let $W$ be the cone in $\PP^4$ over the curve
$C$ with vertex $o_1$. Then we get an absurd inequality
\[\ell^3=Q_1\cdot Q_2\cdot W\geq |\Lambda|>\ell^3,\] where each $Q_i$ is a
general member of $\MMM$. Therefore, at most $k(d-1)$ points of
$\alpha(\Sigma)$ can lie on a  curve of degree $k$ in $\PP^3$ if
$d\geq k^2+1$.

For the projection  $\beta :\PP^3\dasharrow \PP^2$ from a generic
point $o_2\in\PP^3$, we apply the proof of
Corollary~\ref{corollary:conic} to $\alpha(\Sigma)$. This
completes our proof because the general projection $\phi_4$ is the
composite of the projections $\alpha$ and $\beta$.
\end{proof}

\subsection{Base-point-freeness}
It is a classical result that if $6$ points $p_1, \cdots,
p_6\in\PP^2$ in   general position are blown up, then the complete
linear system on the blow-up corresponding to
$|\mathcal{O}_{\PP^2}(3)-p_1-\ldots-p_6|$ is very ample as well as
base-point-free. This is a key observation to classify del Pezzo
surfaces.

Bese's paper \cite{Bes83} developed this observation to points on
$\PP^2$ in less general position and various divisors. The result
however turned out to have a considerable generalization. E.~Davis
and A.~Geramita obtained a very ampleness and a
base-point-freeness theorems on blow-ups of $\PP^2$ via the
ideal-theoretic route that are more powerful than Bese's.

The theorem below is a special case of the paper \cite{DaGe88}
which  provides a strong enough tool for us to study the
base-point-freeness of linear systems of certain types on blow-ups
of $\PP^2$.

\begin{theorem}
\label{theorem:Bese} Let $\pi: Y\to\PP^2$ be the blow up at points
$p_1, \cdots, p_s$ on $\PP^2$. Then the linear system
$|\pi^*(\mathcal{O}_{\PP^2}(d))-\sum^{s}_{i=1}E_i|$ is
base-point-free for all $s\leq \max\{m(d+3-m)-1, m^2\}$, where
$E_i=\pi^{-1}(p_i)$, $d\geq 3$, and $m=\lfloor
\frac{d+3}{2}\rfloor$, if the set
$\Gamma=\{p_1, \cdots, p_s\}$ satisfies the following:\\
\[\mbox{no $k(d+3-k)-1$ points of $\Gamma$ lie on a single curve of degree
$k$, $1\leq k\leq m$.}\]
\end{theorem}

In the case $d=3$ Theorem~\ref{theorem:Bese} is the well known
result on the base-point-freeness of the anticanonical linear
system of a weak del Pezzo surface of degree $9-s\geqslant 2$. The
theorem above  immediately implies the following:

\begin{corollary}
\label{corollary:Bese} Let $\Gamma=\{p_1, \cdots, p_s\}$ be a
finite set of points in $\mathbb{P}^{2}$. For a given positive
integer $d\geq 3$, if $s\leq \max\{m(d+3-m)-1, m^2\}$ and no
$k(d+3-k)-1$ points of the set $\Gamma$ lie on a curve  of degree
$k\leqslant m$ in $\PP^2$, where $m=\lfloor \frac{d+3}{2}\rfloor$,
then for each point $p_i\in \Gamma$
 there is a curve in $\mathbb{P}^{2}$ of degree $d$ that contains
all the points of the set $\Gamma$ except the point $p_i$.
\end{corollary}

The corollary above is the most important tool for this paper.
Also, it makes us propose Conjecture~\ref{conjecture:projections}.

\subsection{Basic properties of nodes}
As explained at the beginning, the ranks of the 4th homology
groups of nodal hypersurfaces in $\PP^4$ are strongly related to
the number of nodes and their position. Even though the number of
nodes are given in our problem, it is necessary that we should
study their position.

\begin{lemma}\label{lemma:plane}
Let $V$ be a nodal hypersurface of degree $n$ in $\PP^4$.
\begin{enumerate}
\item  A curve of degree $k$ in $\PP^4$ contains at most $k(n-1)$
nodes of $V$.

\item If a $2$-plane contains $\frac{n(n-1)}{2}+1$ nodes of $V$,
then the $2$-plane  is contained in $V$.
\end{enumerate}
\end{lemma}

\begin{proof} Suppose that the hypersurface $V$ is defined by an equation
$F(x_0, x_1, x_2, x_3,x_4)=0$. Then the singular locus of $V$ is
contained in a generic hypersurface $V'=(\Sigma\lambda_i
\frac{\partial F}{\partial x_i}=0)$ of degree $n-1$. Let $C$ be a
curve of degree $k$ in $\PP^4$. Since the hypersurface $V$ has
only isolated singularities, the curve $C$ cannot be contained in
$V'$. Because the intersection number of  the hypersurface $V'$
and the curve $C$ is $k(n-1)$, the curve $C$ contains at most
$k(n-1)$ singular points of $V$.

Let $\Pi$ be a $2$-plane not contained in $V$. Then $V'\cap \Pi$
is a curve of degree $(n-1)$ not contained in $V$.  Therefore, the
curve $V'\cap \Pi$ cannot meet $V$ at more than $\frac{n(n-1)}{2}$
nodes of $V$.
\end{proof}

\begin{lemma}\label{lemma:plane-number-of-nodes}
If a nodal hypersurface $V$ of degree $n$ in $\PP^4$ contains a
$2$-plane, then it has at least $(n-1)^2$ singular points.
\end{lemma}

\begin{proof} Suppose that $V$ contains a $2$-plane $\Pi$. Then take two
general hyperplanes  $H_1$ and $H_2$ passing through the $2$-plane
$\Pi$. Then for each $i=1, 2$ we have  the residual surface $Q_i$
of degree $n-1$ in $H_i$ such that $H_i\cap V=\Pi\cup Q_i$. Then
the $(n-1)^2$ intersection points of $Q_1\cap Q_2\cap \Pi$ are
singular points of $V$. \end{proof}

\subsection{Simple tools}
To prove Theorem~\ref{theorem:main}, we need various tools to
handle hypersurfaces of  certain degrees and a finite number of
points in $\PP^4$.
\begin{lemma}\label{lemma-of-numbers-lines}
Let $\Sigma=\{p_1, \cdots, p_r\}$ be a set of $r$ points in
$\PP^N$, $N\geq 2$. Let $p$ be a point in $\PP^N\setminus\Sigma$.
Suppose that no $m+1$ points of $\Sigma$ lie on a single line with
the point $p$. Then there are at least $\min\{r-m, \lfloor
\frac{r}{2} \rfloor\}$ mutually disjoint pairs of points in
$\Sigma$ such that each pair determines a line not containing the
point $p$.

\end{lemma}

\begin{proof} We may assume that the points $p_1, \cdots, p_m$, and $p$
are on a single line $L$.

First, we suppose  $m\geq r-m$. We then obtain $r-m$ such pairs by
choosing one point from $\Sigma\cap L$ and the other from
$\Sigma\setminus L$. Obviously, these pairs determine lines not
passing through the point $p$.

Next, we suppose $m<r-m$. We can then find $\lfloor
\frac{r-2m}{2}\rfloor$ such pairs by choosing two points of
$\Sigma\setminus L$; otherwise $m+1$ points of $\Sigma$ would lie
on a single line. By choosing one point from the remaining points
in $\Sigma\setminus L$ and the other from $L$ we also obtain $m$
such pairs.  The number of the pairs that we obtain is $m +\lfloor
\frac{r-2m}{2}\rfloor=\lfloor \frac{r}{2}\rfloor$. \end{proof}

\begin{lemma}\label{lemma:two-point-cone}
Let $H$ be a hyperplane in $\PP^N$, $N\geq 3$ and $X$ be a
hypersurface of degree $d\geq 2$ in $H$. Suppose that the
hypersurface $X$ does not contain a point $o\in H$. For given two
points $p, q$ in $\PP^N\setminus H$, there is a hypersurface of
degree $d$ in $\PP^N$ such that contains $X$ and two points $p$
and $q$ but not the point $o$.
\end{lemma}

\begin{proof}
Take a generic $2$-plane $\Pi$ passing through the points $p$ and
$q$ in $\PP^N$. Then $\Pi$ meets $X$ at at least two points, say
$p'$ and $q'$. Then the line determined by  $p$ and $p'$ and the
line determined by  $q$ and $q'$ meet at a point $v$. Then the
cone over $X$ with vertex $v$ is a hypersurface of degree $d$ in
$\PP^N$ containing $X$ and two points $p$ and $q$ but does not
contain the point $o$.
\end{proof}

\begin{lemma}\label{lemma:sweeping-method}
Let $\Lambda$ and $\Delta$ be  disjoint finite sets of points in
$\PP^N$, $N\geq 3$ and let $p$ be a point in
$\PP^N\setminus(\Lambda\cup\Delta)$. Suppose that $D_0$ be a
hypersurface of degree $k$ in $\PP^N$ containing the set $\Lambda$
but not the point $p$. If for each point $q\in\Delta$ there is a
hypersurface $D_q$ of degree $k$ in $\PP^N$ containing the set
$(\Lambda\cup\Delta\cup\{p\})\setminus\{q\}$ but not the point
$q$, then there is a hypersurface of degree $k$ such that passes
through $\Lambda\cup\Delta$ but not the point $p$ in $\PP^N$.
\end{lemma}
\begin{proof}
Suppose that the hypersurface $D_0$ is defined by a homogeneous
polynomial $g(x_0, \cdots, x_N)$. Also, we suppose that each
hypersurface $D_q$ is defined by a homogeneous polynomial
$f_q(x_0, \cdots, x_N)$. Then $g(p)\ne 0$ and $f_q(p)=0$ for each
$q\in\Delta$. Furthermore, $f_q(q')\ne 0$ for some $q'\in\Delta$
if and only if $q=q'$. There is a complex number $c_q$ for each
$q\in\Delta$ such that $g(q)+c_qf_q(q)=0$ because $f_q(q)\ne 0$.
Then the hypersurface defined by
\[g(x_0,\cdots, x_N)+\sum_{q\in\Delta}c_qf_q(x_0,\cdots, x_N)=0\]
contains the set $\Lambda\cup\Delta$ but not the point $p$.
\end{proof}

\begin{corollary}\label{corollary:non-vanishing-plus-two}
Let $\MMM$ be a linear system consisting of hypersurfaces of
degree $k\geq 2$ on $\PP^N$, $N\geq 3$. If the base locus
$\Lambda$ of the linear system $\MMM$ is zero-dimensional, then
for two distinct points $p, q$ in $\PP^N\setminus\Lambda$ and a
point $o$ in $\Lambda$, there is a hypersurface of degree $N(k-1)$
such that passes through $\Lambda\cup\{p,q\}\setminus \{o\}$ but
not the point $o$ in $\PP^N$.
\end{corollary}
\begin{proof}
By Lemma~\ref{lemma:non-vanishing}, there is a hypersurface $D_0$
of degree $N(k-1)$ in $\PP^N$ that passes through all the points
of $\Lambda$ except the point $o$.  Let $D$ be a general member in
$\MMM$.  We  choose a hyperplane $H_p$ in $\PP^N$ that passes
through the point $p$ but not the point $q$. We  also choose a
hyperplane $H_q$ that passes through the point $q$ but not the
point $p$. We  then apply Lemma~\ref{lemma:sweeping-method} to the
hypersurfaces $D_0$, $D+(N(k-1)-k)H_p$, and $D+(N(k-1)-k)H_q$ .
\end{proof}

The result below is originally due to J. Edmonds (\cite{Ed65}). It
can help us to make our proofs  simpler.
\begin{theorem}\label{Eisenbud-Koh}
Let $\Sigma$ be a set of points in $\mathbb{P}^N$ and let $d \geq
2$ be an integer. If no $dk+2$  points of $\Sigma$ lie in a
$k$-plane for all $k \geq 1$, then the set $\Sigma$ is $d$-normal
in $\PP^N$.
\end{theorem}

\begin{proof} See \cite{EiJ87}. \end{proof}

\section{Conjectural proof}
\label{section:hypersurfaces}

In this section, we prove Conjecture~\ref{conjecture:factoriality}
under the assumption that Conjectures~\ref{conjecture:projections}
is true.

Let $V_n$ be a nodal hypersurface of degree $n$ in $\PP^4$.
Suppose that $|\Sing(V_n)|<(n-1)^{2}$ and $n\geqslant 4$. Fix a
point $p\in\Sing(V_n)$. And then put
$\Gamma=\Sing(V_n)\setminus\{p\}$. To prove the factoriality of
$V_n$ it is enough to construct a hypersurface of degree $2n-5$
that contains all the points of the set $\Gamma$ and does not
contain the point $p$.

We suppose that Conjecture~\ref{conjecture:projections} holds. Let
$\phi_4:\PP^4\dasharrow \Pi$ be the projection from a general line
$L$ in $\PP^4$, where $\Pi$ is a $2$-plane in $\PP^4$. It then
follows from  Conjecture~\ref{conjecture:projections} and
Lemma~\ref{lemma:plane} that the set $\phi_4(\Gamma)$ satisfies
the condition
\begin{equation}\label{equation:hypersurface-condition}\mbox{\emph{no $k(n-1)+1$ points of $\phi_4(\Gamma)$
lie on a curve of degree $k$ on $\Pi$ for each $k\geq 1$.}}
\end{equation}

\begin{remark}\label{remark:line-conic2}
It follows from Corollaries~\ref{corollary:line} and
\ref{corollary:conic2} that the
condition~(\ref{equation:hypersurface-condition}) holds for
$k\leq\sqrt{n-1}$ without Conjecture~\ref{conjecture:projections}.
\end{remark}

\begin{lemma}\label{lemma-of-after-projection2}
For each $1\leq k\leq n-1$, any curve of degree $k$ on $\Pi$
cannot contain $k(2n-2-k)-1$ points of $\phi_4(\Gamma)$.
\end{lemma}
\begin{proof} It is easy to check that $k(n-1)\leq k(2n-2-k)-2$ if $k<n$.
\end{proof}

\begin{lemma}\label{lemma-final-step2}
There is a curve of degree $2n-5$ on $\Pi$ which passes through
all the points of $\phi_4(\Gamma)$ but not the point $\phi_4(p)$.
\end{lemma}
\begin{proof} It immediately follows from Corollary~\ref{corollary:Bese} and
Lemma~\ref{lemma-of-after-projection2}.
\end{proof}

\begin{proposition}\label{proposition:conjectural2}
Conjecture~\ref{conjecture:projections} implies
Conjecture~\ref{conjecture:factoriality}.
\end{proposition}

\begin{proof}
 Lemma~\ref{lemma-final-step2} implies that there is a curve $C$
of degree $2n-5$ on $\Pi$ which passes through all the points of
$\phi_4(\Gamma)$ but not the point $\phi_4(p)$. We take the cone
over $C$ with vertex $L$. The cone then contains all the points of
$\Gamma$ but not the point $p$. It implies that if the
hypersurface $V_n$ has  $s<(n-1)^2$ singular points, then these
$s$ points impose $s$ linearly independent conditions on
homogeneous forms of degree $2n-5$ on $\PP^4$. Consequently, the
rank of $4$th singular homology group of $V_n$ is $1$ by
Theorem~\ref{theorem-of-Cynk}, which completes the proof.
\end{proof}

\section{Proof of Theorem~\ref{theorem:main}}
In this section, we prove Theorem~\ref{theorem:main}. Let $V_n$ be
a nodal hypersurface of degree $n$ in $\PP^4$ with at most
$(n-1)^2-1$ nodes. However, as we will see in the proofs, we may
assume that the hypersurface $V_n$ has exactly $(n-1)^2-1$ nodes.
To prove the factoriality of the hypersurface $V_n$, for an
arbitrary point $p\in\Sing(V_n)$, we have to construct a
hypersurface of degree $(2n-5)$ in $\PP^4$ that contains the set
$\Sing(V_n)$ except  the point $p$.

\subsection{Quintic hypersurfaces} \label{section:quintics}
Let $V_5$ be a nodal quintic hypersurface in $\PP^4$ with $15$
nodes. A line can contain at most $4$ nodes by
Lemma~\ref{lemma:plane}. If a $2$-plane $\Pi$ contains $12$ nodes
of $V_5$, then $\Pi$ is contained in $V_5$ by
Lemma~\ref{lemma:plane}. It then follows from
Lemma~\ref{lemma:plane-number-of-nodes} that $V_5$ must have at
least $16$ nodes, which contradicts our assumption. Therefore, a
$2$-plane can contain at most $11$ nodes of $V_5$. Therefore, the
set of nodes of $V_5$ satisfies the condition for $d=5$ in
Theorem~\ref{Eisenbud-Koh} and hence the nodal quintic
hypersurface $V_5$ is factorial.

\subsection{Sextic hypersurfaces}

Let $V_6$ be a nodal sextic hypersurface in $\PP^4$ with $24$
nodes. We denote the set of  nodes of $V_6$ by $\Sigma$. If a
$2$-plane $\Pi$ contains $16$ nodes of $V_6$, then $\Pi$ is
contained in $V_6$ by Lemma~\ref{lemma:plane}. It then follows
from Lemma~\ref{lemma:plane-number-of-nodes} that $V_6$ must have
at least $25$ nodes, which contradicts our assumption. Therefore,
a $2$-plane  contains at most $15$ nodes of $V_6$.

\begin{proposition}\label{proposition:23points}
If no $23$ points of $\Sigma$ lie on a single $3$-plane, then the
hypersurface $V_6$ is factorial.
\end{proposition}

\begin{proof} Since no $6$ points of $\Sigma$ lie on a single line and no
$16$ points of $\Sigma$ lie on a single $2$-plane, the set
$\Sigma$ satisfies the condition for $d=7$ of
Theorem~\ref{Eisenbud-Koh}. Therefore, the set $\Sigma$ is
$7$-normal in $\PP^4$ and hence the hypersurface $V_6$ is
factorial.
\end{proof}

Pick an arbitrary point $p$ in $\Sigma$ and then we denote the set
$\Sigma\setminus\{p\}$ by $\Gamma$. To prove the factoriality of
$V_6$, we must find a hypersurface of degree $7$ in $\PP^4$ that
contains the set $\Gamma$ but not the point $p$.

Due to Proposition~\ref{proposition:23points}, we may assume that
at least $23$ points of $\Sigma$ lie in a single $3$-plane $H$.
Furthermore, we may assume that all the $24$ points of $\Sigma$
lie in the $3$-plane $H$ because in what follows we will show that
there is a septic hypersurface in $H$, not in $\PP^4$, that
contains $\Gamma\cap H$ but not the point $p$.

We  consider the projection $\phi_3:H\dasharrow\Pi$ from a generic
point $o$ in $H$, where $\Pi$ is a generic hyperplane of $H$. At
most $5$ points of $\Sigma$ can lie on a single line in $H$ and at
most $10$ points of $\Sigma$ can lie on a conic on $H$.

\begin{lemma}
If there is a set $\Lambda$ of at least $20$ points of $\Gamma$
such that $\phi_3(\Lambda)$ is contained in a cubic curve  $C$ on
$\Pi$, then there is a septic hypersurface in $H$ that contains
the set $\Gamma$ but not the point $p$.
\end{lemma}
\begin{proof}
We may assume that the cubic curve $C$ contains the point
$\phi_3(p)$. If not, then we can easily construct a septic surface
in $H$ that contains $\Gamma$ but not the point $p$. The curve $C$
must be irreducible because a line (a conic, resp.) contains at
most $5$ ($10$, resp.) points of $\phi_3(\Lambda)$ by
Corollaries~\ref{corollary:line} and \ref{corollary:conic}. It
then follows from  Lemma~\ref{lemma:zero-dimensional} that the
linear system of cubic surfaces in $H$ passing through
$\Lambda\cup\{p\}$ has zero-dimensional base locus. Therefore,
applying Corollary~\ref{corollary:non-vanishing-plus-two}, we
obtain a sextic surface $F$ that passes through $22$ points of
$\Sigma$ but not the point $p$. Note that $|\Gamma\setminus F|\leq
1$. By taking a general hyperplane passing through the point in
$\Gamma\setminus F$, we can construct a septic surface in $H$ that
contains the set $\Gamma$ but not the point $p$.
\end{proof}

From now, we suppose that no $20$ points of $\phi_3(\Lambda)$ is
contained in a cubic curve on $\Pi$. And then let us apply the
similar technique as Lemma~\ref{lemma:non-vanishing}, which has
evolved from the papers \cite{Ch04e} and \cite{PaWoo05}, to the
following case.

\begin{lemma}
If the set $\phi_3(\Lambda)$ is contained in a quartic curve  $C$
on $\Pi$, then there is a septic hypersurface in $H$ that contains
the set $\Gamma$ but not the point $p$.
\end{lemma}
\begin{proof}
The curve $C$ must be  irreducible because of our assumption.
Also, we may assume that it contains the point $\phi_3(p)$ as
well. Then the linear system $\MMM$ of quartic hypersurfaces in
$H$ passing through $\Sigma$ has zero-dimensional base locus.
Meanwhile, we have the sextic surface $Y=H\cap V_6$ contains all
the nodes of $V_6$. It may have non-isolated singularities.
However, it is irreducible and reduced; otherwise the hypersurface
$V_6$ would have more than $24$ nodes. Choose a general enough
surface $S'$ in $\MMM$. Then it is smooth in the outside of the
base locus of $\MMM$ and hence it is normal. Also, the surface $Y$
gives us a reduced divisor $D_6\in|\OOO_{S'}(6)|$ on $S'$. Let
$D_4$ be a divisor in $|\OOO_{S'}(4)|$ given by a general member
of $\MMM$. We then consider the $\QQ$-divisor
$D=(1-\epsilon)D_6+2\epsilon D_4$, where $\epsilon$ is
sufficiently small enough rational number. Then it is easy to
check that the support of $\LLL(S',D)$ is zero-dimensional and
contains $\Sigma$. Use Theorem~\ref{theorem:Shokurov} to obtain
$H^1(S',\III(S',D)\otimes\OOO_{S'}(7))=0$. Therefore, there is a
divisor in $|\OOO_{S'}(7)|$ that contains $\Gamma$ but not the
point $p$. The exact sequence \[0\to H^0(\PP^3,
\OOO_{\PP^3}(3))\to H^0(\PP^3, \OOO_{\PP^3}(7))\to H^0(S',
\OOO_{S'}(7))\to 0\] completes the proof.
\end{proof}

\begin{lemma}
If the set $\phi_3(\Gamma)$ is not contained in any quartic curve
$C$ on $\Pi$, then there is a septic hypersurface in $H$ that
contains the set $\Gamma$ but not the point $p$.
\end{lemma}
\begin{proof}
In this case, the set $\phi_3(\Lambda)$ satisfies the condition
for $d=7$ in Theorem~\ref{corollary:Bese}. Therefore, there is a
septic curve on $\Pi$ that contains the set $\phi_3(\Lambda)$ but
not the point $\phi_3(p)$. Then the cone over the septic curve
with vertex $o$ is a septic hypersurface in $H$ that contains
$\Gamma$ but not the point $p$.
\end{proof}

Consequently, for an arbitrary point $p\in\Sigma$, we can find a
septic hypersurface in $\PP^4$ that contains $\Gamma$ but not the
point $p$. Therefore, the rank of the 4th integral homology group
of $V_6$ is $1$ and hence $V_6$ is factorial.
\subsection{Septic hypersurfaces}

Let $V_7$ be a nodal septic hypersurface in $\PP^4$ with $35$
nodes. We again denote the set of  nodes of $V_7$ by $\Sigma$. Fix
a point $p$ in $\Sigma$. We denote the set $\Sigma\setminus\{p\}$
by $\Gamma$. To prove the factoriality of $V_7$, we have to
construct a hypersurface of degree $9$ in $\PP^4$ that contains
the set $\Gamma$ but not the point $p$.

First of all, it follows from
Lemma~\ref{lemma:plane-number-of-nodes} that a $2$-plane contains
at most $21$ points of $\Sigma$. Suppose that  there is a
$2$-plane $\Pi$ containing at least $20$ points of $\Sigma$. The
$2$-plane $\Pi$ is not contained in $V_7$; otherwise the
hypersurface $V_7$ would have at least $36$ nodes. We then
consider the projection $\phi_4:\PP^4\dasharrow \Pi$ from a
generic line $L$. By Corollaries~\ref{corollary:line}
and~\ref{corollary:conic2} a line on $\Pi$ contains at most $6$
points of $\phi_4(\Sigma)$ and a conic on $\Pi$ contains at most
$12$ points of $\phi_4(\Sigma)$.

\begin{proposition}
If there is a $2$-plane $\Pi$ containing at least $20$ points of
$\Sigma$, then for the projection $\phi_4:\PP^4\dasharrow \Pi$
from a generic line $L$ the set $\phi_4(\Gamma)$ satisfies the
following:
\begin{enumerate}
\item A line on $\Pi$ contains at most $6$ points of
$\phi_4(\Gamma)$.

\item A conic on $\Pi$ contains at most $12$ points of
$\phi_4(\Gamma)$.

\item A cubic on $\Pi$ contains at most $25$ points of
$\phi_4(\Gamma)$.

\item A quartic on $\Pi$ contains at most $30$ points of
$\phi_4(\Gamma)$.

\item A quintic on $\Pi$ contains at most $33$ points of
$\phi_4(\Gamma)$.

\item A sextic on $\Pi$ contains at most $34$ points of
$\phi_4(\Gamma)$.
\end{enumerate}
\end{proposition}

\begin{proof} The first and the second statements follow from
Corollaries~\ref{corollary:line} and~\ref{corollary:conic2}. And
the last statement is obvious because $|\phi_4(\Gamma)|=34$.

For a cubic, we suppose that there is a cubic $C$ on $\Pi$ that
contains $26$ points $\phi_4(p_1),\cdots,\phi_4(p_{26})$ of
$\phi_4(\Gamma)$. The cubic $C$ must be irreducible because of the
first and the second statements. It then follows from
Lemma~\ref{lemma:zero-dimensional} that the base locus of the
linear system $\MMM$ of cubic hypersurfaces in $\PP^4$ containing
the points $p_1,\cdots, p_{26}$ is zero-dimensional and hence the
restricted linear system $\MMM|_\Pi$ of the linear system $\MMM$
to the $2$-plane $\Pi$ also has zero-dimensional base locus. Since
we have at most $15$ points of $\Sigma$ in the outside of $\Pi$,
at least $11$ points of $p_1,\cdots, p_{26}$ belong to $\Pi$.
Therefore, there is an irreducible  cubic curve $D$ in $\MMM|_\Pi$
that is not contained in $V_7$ but passing through $11$ nodes of
$V_7$. However, this is impossible because $21=D\cdot V_7\geq
11\cdot 2=22$.

If we have a quartic in $\Pi$ containing $31$ points of
$\phi_4(\Gamma)$, then in the same way, we can find an irreducible
quartic curve not contained in $V_7$ but passing through $16$
nodes of $V_7$, which is absurd.

Finally, if there is a quintic in $\Pi$ containing $34$ points of
$\phi_4(\Gamma)$, then in the same way we can find an irreducible
quintic curve not contained in $V_7$ but passing through $19$
nodes of $V_7$, which is also impossible. \end{proof}

\begin{corollary}
If there is a $2$-plane containing at least $20$ points of
$\Sigma$, the hypersurface $V_7$ is factorial.
\end{corollary}

\begin{proof} The proposition above shows the set
$\phi_4(\Gamma)$ satisfies the condition for $d=9$ in
Corollary~\ref{corollary:Bese} and hence there is a curve $C$ of
degree $9$ on $\Pi$ passing through all the point of
$\phi_4(\Gamma)$ but not the point $\phi_4(p)$. The cone over the
curve $C$ with vertex $L$ shows that the set $\Sigma$ is
 $9$-normal in
$\PP^4$. Therefore, the hypersurface $V_7$ is factorial.
\end{proof}

From now on, we suppose that a $2$-plane contains at most $19$
points of $\Sigma$. If a hyperplane in $\PP^4$ contains at most
$28$ points of $\Sigma$, then the set $\Sigma$ satisfies the
condition for $d=9$ in Theorem~\ref{Eisenbud-Koh} and hence it is
$9$-normal and $V_7$ is factorial. Therefore, we suppose that a
hyperplane $H$ in $\PP^4$ contains at least $29$ points of
$\Sigma$. And let $\Sigma'=\Sigma\cap H$ and
$\Sigma''=\Sigma\setminus H$. We always assume that the point $p$
is contained in $\Sigma'$ because, if not, then we can easily
construct a hypersurface of degree $9$ in $\PP^4$ containing
$\Sigma$ except the point $p$.

We consider the projection $\alpha :H\dasharrow \Pi$ from a
generic point $o_1\in H$, where $\Pi$ is a general $2$-plane in
$H$.

\begin{lemma}
If there is a set $\Lambda$ of at least $26$ points of
$\Sigma'\setminus\{p\}$ such that $\alpha(\Lambda)$ is contained
in a cubic curve $C$ on $\Pi$, then there is a hypersurface of
degree $9$ that contains the set $\Gamma$ but not the point $p$.
\end{lemma}
\begin{proof}
Note that the curve $C$ is irreducible and
$m=|\Gamma\setminus\Lambda|\leq 8$.

 Suppose that the curve $C$
does not contain the point $\alpha(p)$. Because a line has at most
$6$ points of $\Sigma$, there is a hypersurface of degree
$\min\{m-5, \lfloor \frac{m}{2}\rfloor\}\leq 4$ in $\PP^4$ that
passes through $\Gamma\setminus\Lambda$ but not the point $p$ by
Lemma~\ref{lemma-of-numbers-lines}. Therefore, we can easily
construct a  hypersurface of degree $9$ that passes through
$\Gamma$ but not  the point $p$.

Suppose that the curve $C$ contains also the point $\alpha(p)$.
Pick two points $p_1$ and $p_2$ from $\Gamma\setminus\Lambda$ in
such a way that $\Gamma\setminus (\Lambda\cup\{ p_1, p_2\})$ is
contained in a cubic hypersurface $F_1$ in $\PP^4$ not containing
the point $p$, which is possible because of
Lemma~\ref{lemma-of-numbers-lines}. The linear system of cubic
hypersurfaces in $H$ containing $\Lambda\cup \{p\}$ has zero
dimensional base locus. Therefore, there is a sextic hypersurface
$F_2$  in $\PP^4$ that passes through $\Lambda$ and the points
$p_1$ and $p_2$ but not the point $p$ by
Corollary~\ref{corollary:non-vanishing-plus-two}. Then the nonic
hypersurface $F=F_1+ F_2$  contains all the point of $\Sigma$
except the point $p$.
\end{proof}

We may assume that no $26$ points of
$\alpha(\Sigma'\setminus\{p\})$ lie on a cubic curve on $\Pi$. If
a cubic curve on $\Pi$ contains more than $18$ points of
$\alpha(\Sigma'\setminus\{p\})$  then it must be irreducible.

\begin{lemma}\label{lemma:cubic-excluding-six-points}
If a cubic curve $C$ on $\Pi$ contains $22$ points of
$\alpha(\Sigma'\setminus\{p\})$, it is unique. \end{lemma}
\begin{proof} Suppose that a cubic curve $C'$ on $\Pi$ contains at
least $22$ points of $\alpha(\Sigma'\setminus\{p\})$, then it
meets $C$ at $22-(34-22)=10$ points and hence $C=C'$.
\end{proof}

To prove the factoriality of $V_7$, we will consider the following
five cases.
\bigskip

\emph{Case 1.} $|\Sigma'|=29$.\bigskip

Because no $22$ points of $\Sigma$ are contained in a $2$-plane,
Lemma~\ref{lemma:cubic-excluding-six-points} enables us to choose
six points $p_1, p_2, \cdots, p_6$ from the set
$\Sigma'\setminus\{p\}$ in such a way that \begin{itemize} \item
no $20$ points of $\alpha(\Sigma'\setminus\{p, p_1,\cdots p_6\})$
lie on a single cubic on $\Pi$; \item $p_1$, $p_2$, and $p_3$ lie
on a $2$-plane $\Pi_1$ not containing the point $p$; \item$p_4$,
$p_5$, and $p_6$ lie on a $2$-plane $\Pi_2$ not containing the
point $p$.
\end{itemize}

Then the set $\alpha(\Sigma'\setminus\{p, p_1,\cdots, p_6\})$
satisfies the condition of Corollary~\ref{corollary:Bese} for
$d=7$. Therefore, there is a septic curve $D$ on $\Pi$ containing
$\alpha(\Sigma'\setminus\{p, p_1,\cdots, p_6\})$ but not the point
$p$. Then the cone over $D$ with vertex $o_1$ is the septic
surface $\bar D$ in $H$ containing $\Sigma'\setminus\{p,
p_1,\cdots p_6\}$ but not the point $p$. Choose two points $q_1$
and $q_2$ from $\Sigma''$. By Lemma~\ref{lemma:two-point-cone}, we
can construct a hypersurface $\tilde D$ in $\PP^4$ containing
$\Sigma'\setminus\{p, p_1,\cdots, p_6\}$ and $q_1$ and $q_2$ but
not the point $p$. Now we choose another two points $q_3$ and
$q_4$ from $\Sigma''$. Let $H_3$ be a hyperplane passing through
the point $q_3$ but not $q_4$ and $H_4$ a hyperplane passing
through the point $q_4$ but not $q_3$. Apply
Lemma~\ref{lemma:sweeping-method} to the hypersurfaces $\tilde D$,
$H+ H'+5H_4$, and $H+ H'+5H_3$, where $H'$ is a general hyperplane
in $\PP^4$ passing through $q_1$ and $q_2$, to  obtain a septic
hypersurface in $\PP^4$ passing through $\Sigma'\setminus\{p,
p_1,\cdots p_6\}$ and $\{q_1, q_2, q_3, q_4\}$ but not the point
$p$. By our construction, two $2$-planes $\Pi_1$ and $\Pi_2$ and
the remaining two points of $\Sigma''$ are contained in a
quadratic hypersurface in $\PP^4$ not passing through the point
$p$.
\bigskip

\emph{Case 2.} $|\Sigma'|=30$.\bigskip

For the same reason as in Case 1, we can choose three points
$p_1$, $p_2$, and $p_3$ from the set $\Sigma'\setminus\{p\}$ in
such a way that
\begin{itemize} \item
no $23$ points of $\alpha(\Sigma'\setminus\{p, p_1,p_2, p_3\})$
lie on a single cubic on $\Pi$; \item $p_1$, $p_2$, and $p_3$ lie
on a $2$-plane $\Pi_1$ not containing the point $p$.
\end{itemize}
Then the set $\alpha(\Sigma'\setminus\{p, p_1,p_2, p_3\})$
satisfies the condition of Corollary~\ref{corollary:Bese} for
$d=8$. Therefore, there is an octic curve $D$ on $\Pi$ containing
$\alpha(\Sigma'\setminus\{p, p_1,p_2, p_3\})$ but not the point
$p$. Then the cone over $D$ with vertex $o_1$ is an octic surface
$\bar D$ in $H$ containing $\Sigma'\setminus\{p, p_1,p_2,p_3\}$
but not the point $p$. Choose two points $q_1$ and $q_2$ from
$\Sigma''$. By Lemma~\ref{lemma:two-point-cone}, we can construct
a hypersurface $\tilde D$ in $\PP^4$ containing
$\Sigma'\setminus\{p, p_1,p_2, p_3\}$ and $q_1$ and $q_2$. Now we
choose another two points $q_3$ and $q_4$ from $\Sigma''$. As in
Case~1, we use Lemma~\ref{lemma:sweeping-method} to get an octic
hypersurface in $\PP^4$ passing through $\Sigma'\setminus\{p,
p_1,p_2, p_3\}$ and $\{q_1, q_2, q_3, q_4\}$ but not the point
$p$. Also, the $2$-plane $\Pi_1$ and the remaining one point
$\Sigma''$ are contained in a hyperplane in $\PP^4$ not passing
through the point $p$.

\bigskip
\emph{Case 3.} $|\Sigma'|=31$. \bigskip

Then the set $\alpha(\Sigma'\setminus\{p\})$ satisfies the
condition of Corollary~\ref{corollary:Bese} for $d=9$. Therefore,
there is a nonic curve $D$ on $\Pi$ containing
$\alpha(\Sigma'\setminus\{p\})$ but not the point $p$. Then the
cone over $D$ with vertex $o_1$ is a nonic surface $\bar D$ in $H$
containing $\Sigma'\setminus\{p\}$ but not the point $p$. Choose
two points $q_1$ and $q_2$ from $\Sigma''$. By
Lemma~\ref{lemma:two-point-cone}, we can construct hypersurface
$\tilde D$ in $\PP^4$ containing $\Sigma'\setminus\{p, p_1,p_2,
p_3\}$ and $q_1$ and $q_2$. Note that $|\Sigma''\setminus\{q_1,
q_2\}|=2$.  As in the previous, we use
Lemma~\ref{lemma:sweeping-method} to construct a nonic
hypersurface in $\PP^4$ passing through $\Gamma$ but not the point
$p$.

\bigskip
\emph{Case 4.} $32\leq |\Sigma'|\leq 34$.\bigskip

Suppose that no $31$ points of $\alpha(\Sigma'\setminus\{p\})$ lie
on a single quartic curve on $\Pi$. The set
$\alpha(\Sigma'\setminus\{p\})$ then satisfies  the condition of
Corollary~\ref{corollary:Bese} for $d=9$. Therefore, as in Case 3,
we can find a nonic hypersurface that we need.

Suppose that there is a set $\Lambda$ of at least $31$ points of
$\Sigma'\setminus\{p\}$ such that a quartic curve $C$ on $\Pi$
contains $\alpha(\Lambda)$. If the curve $C$ does not contain the
point $\alpha(p)$, then we can easily construct a nonic
hypersurface in $\PP^4$ that we are looking for. Therefore,
 we may assume that the point $\alpha(p)$ also belongs to $C$.
 Then it follows from Lemmas~\ref{lemma:non-vanishing}
 and~\ref{lemma:zero-dimensional} that there is a nonic hypersurface
 in $H$ passing through $\Lambda$ but not the point $p$. Then,
 using Lemmas~\ref{lemma:two-point-cone} and \ref{lemma:sweeping-method},
 we can construct a nonic hypersurface in $\PP^4$ that passes
 through all the point of $\Gamma$ but not the point
 $p$.

\bigskip
\emph{Case 5.} $|\Sigma'|=35$.

Suppose that there is a set $\Lambda$ of at least $31$ points of
$\Gamma=\Sigma'\setminus\{p\}$ such that $\alpha(\Lambda)$ is
contained in a quartic curve $C$ on $\Pi$. The curve $C$ must be
irreducible.  We also assume that the curve $C$ contains the point
$\alpha(p)$. Then the base locus of the linear system $\MMM$ of
quartic surfaces on $H$ containing $\Lambda\cup\{p\}$ is
zero-dimensional by Lemma~\ref{lemma:zero-dimensional}.  Let $B$
be the support of the base locus of the linear system $\MMM$ and
$\bar\Sigma=\Sigma\setminus B$. Note that $\Lambda\cup\{p\}\subset
B$. It follows from Lemma~\ref{lemma:non-vanishing} that the set
$B$ is $9$-normal in $H$. There is a nonic hypersurface $F$ in $H$
that contains $B\setminus\{p\}$ but not the point $p$. Because
$|\bar\Sigma|\leq 3$, for each $q\in \bar\Sigma$ there is a
quintic hypersurface $Q_q$ in $H$ such that contains the set
$\bar\Sigma$ but not the point $q$.
 Choose a general element $Q$ from the linear system $\MMM$.
We then apply Lemma~\ref{lemma:sweeping-method} to the nonic
hypersurfaces $F$ and $Q+Q_q$ to obtain a nonic hypersurface
passing through $\Sigma$ except the point $p$. Therefore, we may
assume that no $31$ points of $\alpha(\Gamma$) lie on a quartic on
$\Pi$.

Unless the set $\alpha(\Gamma)$ lie on a quintic curve on $\Pi$,
we can use Corollary~\ref{corollary:Bese} to get a nonic curve on
$\Pi$ containing the set $\alpha(\Gamma)$ but not the point
$\alpha(p)$, which gives us a nonic hypersurface in $\PP^4$ that
we need.

Finally, we suppose that there is a quintic curve $C_5$ on $\Pi$
that contains $\alpha(\Gamma)$.

The curve $C_5$ is irreducible. Also, we may assume that it
contains the point $\alpha(p)$ as well. Then the linear system
$\DDD$ of quintic hypersurfaces in $H$ passing through $\Sigma$
has zero-dimensional base locus. Meanwhile, we have the septic
surface $Y=H\cap V_7$ contains all the nodes of $V_7$, which  may
have non-isolated singularities. However, it is irreducible and
reduced; otherwise the hypersurface would have more than $35$
nodes. Choose a general enough surface $S'$ in $\DDD$. Then it is
smooth in the outside of the base locus of $\DDD$ and hence it is
normal. Also, the surface $Y$ gives us a reduced divisor
$D_7\in|\OOO_{S'}(7)|$ on $S'$. Let $D_5$ be a divisor in
$|\OOO_{S'}(5)|$ given by a general member of $\DDD$. We then
consider the $\QQ$-divisor $D=(1-\epsilon)D_7+2\epsilon D_5$,
where $\epsilon$ is sufficiently small enough rational number.
Then it is easy to check that the support of $\LLL(S',D)$ is
zero-dimensional and contains $\Sigma$. Using
Theorem~\ref{theorem:Shokurov},  we obtain
$H^1(S',\III(S',D)\otimes\OOO_{S'}(9))=0$. Therefore, there is a
divisor in $|\OOO_{S'}(9)|$ that contains $\Gamma$ but not the
point $p$. Because the sequence \[0\to H^0(\PP^3,
\OOO_{\PP^3}(4))\to H^0(\PP^3, \OOO_{\PP^3}(9))\to H^0(S',
\OOO_{S'}(9))\to 0\] exact, we are done.

\end{document}